\documentclass{amsart}
\usepackage{amsmath,amssymb,color}
\usepackage[all]{xy}

\newcommand{\Bh}{{\mathcal B}}
\newcommand{\Ch}{{\mathcal C}}

\newcommand{\Hh}{{\mathcal H}}

\newcommand{\Mh}{{\mathcal M}}
\newcommand{\Nh}{{\mathcal N}}

\newcommand{\be}{\mathbf{1}}

\newcommand{\halb}{\frac{1}{2}}

\newcommand{\id}{\mathrm{id}}

\newcounter{number}[section]

\newenvironment{nummer}{\refstepcounter{number}
{\noindent\arabic{section}.\arabic{number}}}{}

\newcommand{\bn}{\begin{nummer} \rm}
\newcommand{\en}{\end{nummer}}

\newenvironment{thms}{\noindent {\sc Theorem:} \it}{}

\newenvironment{props}{\noindent {\sc Proposition:} \it}{}
\newenvironment{dfs}{\noindent {\sc Definition:} \it}{}
\newenvironment{cors}{\noindent {\sc Corollary:} \it}{}

\newenvironment{rems}{\noindent {\sc Remark:}}{}

\newenvironment{exss}{\noindent {\sc Examples:} }{}

\newenvironment{nots}{\noindent {\sc Notation:} }{}

\newenvironment{nproof}{\noindent {\sc Proof:}}{\mbox{}\hfill
\rule[-.2ex]{.25em}{1.8ex}}


\begin{document}
\author{Wilhelm Winter}
\author{Joachim Zacharias}

\address{School of Mathematical Sciences,  University of Nottingham, Nottingham, NG7 2RD}
\email{wilhelm.winter@nottingham.ac.uk}

\email{joachim.zacharias@nottingham.ac.uk}


\title[Completely positive maps of order zero]{Completely positive maps of order zero}

\begin{abstract}
We say a completely positive contractive map between two $C^{*}$-algebras has order zero, if it sends orthogonal elements to orthogonal elements. We prove a structure theorem for such maps. As a consequence, order zero maps are in one-to-one correspondence with $*$-homomorphisms from the cone over the domain into the target algebra.  Moreover, we conclude that tensor products of order zero maps are again order zero,  that the composition of an order zero map with a tracial functional is again a tracial functional, and that order zero maps respect the Cuntz relation, hence induce ordered semigroup morphisms between Cuntz semigroups.  
\end{abstract}

\maketitle

\setcounter{section}{-1}

\section{Introduction}   
\label{introduction}

There are various types of interesting maps between $C^{*}$-algebras, all of which can serve as morphisms of a category with objects (a subclass of) the class of all $C^{*}$-algebras. As a first choice, continuous $*$-homomorphisms come to mind, and it follows from spectral theory that in fact any $*$-homomorphism between $C^{*}$-algebras is automatically continuous. At the opposite end of the scale, one might simply consider (bounded) linear maps. It is then a natural question which classes of morphisms interpolate in a reasonable way between linear maps and $*$-homomorphisms. For example, one might ask a linear map to preserve the involution, or even the order structure, i.e., to be self-adjoint or positive, respectively.  In noncommutative topology, it is also often desirable to consider maps which have well-behaved amplifications to matrix algebras; this leads to the strictly smaller classes of completely bounded, or completely positive (c.p.) maps, for example. In contrast, amplifications of $*$-homomorphisms automatically are $*$-homomorphisms. 

Emphasizing the $*$-algebra structure rather than the order structure of a $C^{*}$-algebra, one might also consider Jordan $*$-homomorphisms (amplifications of which again are Jordan $*$-homomorphisms). Another concept, which has recently turned out to be highly useful, but has received less attention in the literature, is that of orthogonality (or disjointness) preserving maps. By this, we mean linear maps which send orthogonal elements to orthogonal elements.   There is a certain degree of freedom here, since one might only ask for orthogonality of supports, or ranges, to be preserved; this distinction becomes irrelevant in the case of c.p.\ maps; we will say an orthogonality preserving c.p.\ map to be of \emph{order zero}. 

In \cite{Wol:disjointness}, Wolff proved a structure theorem for  bounded, linear, self-adjoint, disjointness preserving maps with unital domains: Any such map  is a compression of a Jordan $*$-homomorphism with a self-adjoint element commuting with its image. 

Later (but independently) the first named author arrived at a very similar result for c.p.\ order zero maps in the case of finite-dimensional domains. Any such map is a compression of a $*$-homomorphism with a positive element commuting with its image. Order zero maps with finite-dimensional domains  have been used in \cite{Winter:cpr1}, \cite{Winter:cpr2}, \cite{KirWinter:dr} and \cite{Kir:CentralSequences} as building blocks of noncommutative partitions of unity to define noncommutative versions of topological covering dimension; see \cite{RordamWinter:Z-revisited} and \cite{Winter:localizingEC} for related applications. They will serve a similar purpose in \cite{WinterZac:dimnuc}. However, also order zero maps with more general domains occur in a natural way. To analyze these it will be crucial to have a structure theorem for general c.p.\ order zero maps at hand. 

In the present paper we use Wolff's result to provide such a generalization, see Theorem~\ref{main-result}. Compared to Wolff's theorem, our result produces a stronger statement from stronger hypotheses; it has the additional benefit that it covers the nonunital situation as well.

We obtain a number of interesting consequences from Theorem~\ref{main-result}. For once, it turns out that completely positive contractive (c.p.c.) order zero maps from $A$ into $B$ are in one-to-one correspondence with $*$-homomorphisms from the cone over $A$ into $B$. This point of view also leads to a notion of positive functional calculus for c.p.\ order zero maps. We then observe that tensor products of c.p.\ order zero maps are again order zero; this holds in particular for amplifications of c.p.\ order zero maps to matrix algebras. Moreover, we show that the composition of a c.p.\ order zero map with a tracial functional again is a tracial functional. Finally, we show that (unlike general c.p.\ maps) order zero maps induce ordered semigroup morphisms between Cuntz semigroups. In fact, this observation is one of our motivations for studying order zero maps, since it shows that they provide a natural framework to study the question when maps at the level of Cuntz semigroups can be lifted to maps between $C^{*}$-algebras. For $K$-theory, this problem has been well-studied; it is of particular importance for the classification program for nuclear $C^{*}$-algebras.  While the Cuntz semigroup in recent years also has turned out to be highly relevant for the classification program (cf.\ \cite{EllToms:regularity}, \cite{BroPerToms:cuntz-semigroup}), at this point not even a bivariant version (resembling Kasparov's $KK$-theory)  has been developed. We are confident that our results can be used to build such a theory; this will be pursued in  subsequent work. 

Our paper is organized as follows. In Section~{\ref{orthogonality}} we recall some facts about orthogonality in $C^{*}$-algebras and introduce the notion of c.p.\ order zero maps. In Section~{\ref{main-result-section}}, we prove a unitization result as well as our structure theorem for such maps. We derive a number of corollaries in Section~{\ref{consequences}}.


\newpage

\section{Orthogonality} 
\label{orthogonality}

In this section we recall some facts about orthogonality in $C^{*}$-algebras, introduce the notion of c.p.\ order zero maps, and recall a result of Wolff.

\bn
\begin{nots}
Let $a,b$ be elements in a $C^{*}$-algebra $A$. We say $a$ and $b$ are orthogonal, $a \perp b$, if $ab=ba=a^{*}b=ab^{*}=0$. 
\end{nots}
\en

\bn
\label{orthogonality-remark}
\begin{rems}
In the situation of the preceding definition, note that $a \perp b$  iff $ a^{*}a \perp b^{*}b$, $ a^{*}a \perp bb^{*}$, $aa^{*} \perp b^{*}b$ and $aa^{*} \perp bb^{*}$. 

Note also that, if $a$ and $b$ are self-adjoint, then $a \perp b$ iff $ab=0$.
\end{rems}
\en

\bn
\begin{dfs}
Let $A$ and $B$ be $C^{*}$-algebras and let $\varphi:A \to B$ be a c.p.\ map. We say $\varphi$ has order zero, if, for $a,b \in A_{+}$,
\[
a \perp b \Rightarrow \varphi(a) \perp \varphi(b).
\]  
\end{dfs}
\en

\bn
\begin{rems}
In the preceding definition, we could as well consider general elements $a,b \in A$; this yields the same definition, since    we assume $\varphi$ to be completely positive.  

To see this, note that if $\varphi$ respects orthogonality of arbitrary elements, it trivially has order zero. Conversely, suppose $\varphi$ has order zero, i.e., respects orthogonality of positive elements, and let $a \perp b \in A$ be arbitrary. Then, $a^{*}a \perp b^{*}b$, $a^{*}a \perp bb^{*}$, $aa^{*} \perp bb^{*}$ and $aa^{*} \perp b^{*}b$. We obtain $\varphi(a^{*}a) \perp \varphi(b^{*}b)$, $\varphi(a^{*}a) \perp \varphi(bb^{*})$, $\varphi(aa^{*}) \perp \varphi(bb^{*})$ and $\varphi(aa^{*}) \perp \varphi(b^{*}b)$. But since $\varphi$ is c.p., we have $0 \le \varphi(a^{*}) \varphi(a) \le \varphi(a^{*}a)$, $0 \le \varphi(a) \varphi(a^{*}) \le \varphi(aa^{*})$ (and similarly for $b$ in place of $a$), which yields that $\varphi(a^{*}) \varphi(a) \perp \varphi(b^{*}) \varphi(b)$, $\varphi(a^{*}) \varphi(a) \perp \varphi(b) \varphi(b^{*})$, $\varphi(a) \varphi(a^{*}) \perp \varphi(b) \varphi(b^{*})$ and $\varphi(a) \varphi(a^{*}) \perp \varphi(b^{*}) \varphi(b)$, since orthogonality is a hereditary property. By Remark~\ref{orthogonality-remark}, this implies that $\varphi(a) \perp \varphi(b)$, whence $\varphi$ respects orthogonality of arbitrary elements.
\end{rems}
\en

\bn
\begin{exss}
Any $*$-homomorphism between $C^{*}$-algebras clearly has order zero. More generally, if $\pi:A \to B$ is a $*$-homomorphism and $h \in B$ is a positive element satisfying $[h,\pi(A)]=0$, then $\varphi(\, . \,):= h \pi(\, .\,)$ defines a c.p.\ order zero map. We will show in Theorem~\ref{main-result} that any c.p.\ order zero map is essentially of this form.
\end{exss}
\en

\bn
\label{Wolffs-result}
In \cite{Wol:disjointness}, Wolff defined a bounded linear map to be disjointness preserving, if it is self-adjoint and sends orthogonal self-adjoint elements to orthogonal self-adjoint elements. For the convenience of the reader, we state below  the  main result of that paper, \cite[Theorem~2.3]{Wol:disjointness}. Recall that a Jordan $*$-homomorphism $\pi:A \to B$ between $C^{*}$-algebras is a linear self-adjoint map preserving the Jordan product $a \cdot b = \frac{1}{2} (ab + ba)$; equivalently, one could ask $\pi$ to preserve squares of positive elements, i.e., $\pi(a^{2})=\pi(a)^{2}$ for all $0 \le a \le \be$.

\begin{thms}
Let $A$ and  $B$ be $C^{*}$-algebras, with $A$ unital, and let $\varphi:A \to B$ be a disjointness preserving map. Set $C:= \overline{\varphi(\be_{A}) \{\varphi(\be_{A})\}'}$. Then, $\varphi(A) \subset C$ and there is a Jordan $*$-homomorphism $\pi:A \to \Mh(C)$ from $A$ into the multiplier algebra of $C$ satisfying
\[
\varphi(a) = \varphi(\be_{A}) \pi(a)
\]
for all $a \in A$.
\end{thms}
\en

\section{The main result}
\label{main-result-section}

Below, we prove a unitization result for c.p.\ order zero maps as well as  our main theorem.

\bn
\begin{nots}
Following standard notation, we will write $A^{+}$ for the 1-point unitization of a $C^{*}$-algebra $A$, i.e., $A^{+} \cong A \oplus \mathbb{C}$ as a vector space. If $\varphi:A \to B$ is a c.p.c.\ map into a unital $C^{*}$-algebra $B$, we write $\varphi^{+}:A^{+} \to B$ for  the uniquely determined unital c.p.\ extension of $\varphi$. Recall that, if $\varphi$ is a $*$-homomorphism, then so is $\varphi^{+}$. We denote by $\Mh(A)$ the multiplier algebra of $A$ and by $A^{**}$ its bidual,  identified with the envelopping von Neumann algebra.
\end{nots}
\en

\bn
\label{order-zero-unitization}
\begin{props}
Let $A$ and $B$ be $C^{*}$-algebras, with $A$ nonunital, and let $\varphi:A \to B$ be a c.p.c.\ order zero map. Set $C:= C^{*}(\varphi(A)) \subset B$. 

Then, $\varphi$ extends uniquely to a c.p.c.\ order zero map $\varphi^{(+)}:A^{+} \to C^{**}$.
\end{props}

\begin{nproof}
We may clearly assume that $C$ acts nondegenerate on a Hilbert space $\Hh$.

Choose an increasing approximate unit $(u_{\lambda})_{\lambda \in \Lambda}$ for $A$ and note that 
\begin{equation}
\label{w4}
g:= \mathrm{s.o.} \lim_{\lambda} \varphi (u_{\lambda}) \in C^{**}
\end{equation}
exists in the bidual of $C$ since the $\varphi(u_{\lambda})$ form a bounded, monotone increasing net in $C^{**}$. Define a linear map
\[
\varphi^{(+)}:A^{+} \to C^{**}
\]
by
\[
\varphi^{(+)}(a) := \varphi(a), \, a \in A
\]
and
\begin{equation}
\label{w3}
\varphi^{(+)}(\be_{A^{+}}):= g.
\end{equation}
Note that $\varphi^{(+)}$ is well defined since $A^{+} \cong A \oplus \mathbb{C}$ as a vector space. 

By Stinespring's Theorem, there are a Hilbert space $\Hh_{1}$, a (nondegenerate) $*$-homomorphism 
\[
\sigma : A \to \Bh(\Hh_{1})
\]  
and an operator $v \in \Bh(\Hh,\Hh_{1})$ such that 
\[
v^{*}v \le \be_{\Hh} \mbox{ and } \varphi(a) = v^{*} \sigma(a) v
\]
for $a \in A$. Note that
\begin{eqnarray*}
g & = & \mathrm{s.o.} \lim_{\lambda} \varphi(u_{\lambda}) \\
& = & \mathrm{s.o.} \lim_{\lambda} v^{*} \sigma(u_{\lambda}) v \\
& = & v^{*}(\mathrm{s.o.} (\lim_{\lambda}  \sigma(u_{\lambda}))) v \\
& = & v^{*} \be_{\Hh_{1}} v \\
& = & v^{*}v, 
\end{eqnarray*}
where for the fourth equality we have used that $\sigma$ is nondegenerate. Let 
\[
\sigma^{+}:A^{+} \to \Bh(\Hh_{1})
\] 
be the unitization of $\sigma$, i.e., 
\[
\sigma^{+}(a + \alpha \cdot \be_{A^{+}}) = \sigma(a) + \alpha \cdot \be_{\Hh_{1}} \mbox{ for } a \in A,\, \alpha \in \mathbb{C};
\]
$\sigma^{+}$ is a $*$-homomorphism. We have, for $a \in A$ and $\alpha \in \mathbb{C}$,
\begin{eqnarray*}
v^{*} \sigma^{+}(a + \alpha \cdot \be_{A^{+}}) v & = & v^{*} \sigma(a) v + \alpha \cdot v^{*}v \\
& = & \varphi(a) + \alpha \cdot g \\
& = & \varphi^{(+)}(a + \alpha \cdot \be_{A^{+}}) ,
\end{eqnarray*} 
so $\varphi^{(+)}$ is c.p.c., being a compression of a $*$-homomorphism. 

We next check that $\varphi^{(+)}$ is again an order zero map.  To this end, let 
\[
a + \alpha \cdot \be_{A^{+}} \mbox{ and } b+ \beta\cdot \be_{A^{+}}
\]
in $(A^{+})_{+}$ be orthogonal elements. Since orthogonality passes to quotients, we see that at least one of $\alpha$ and $\beta$ has to be zero. So let us assume $\beta=0$ and note that this implies $b \ge 0$; note also that $a = a^{*}$ and $\alpha \ge 0$.

We have, for each $\lambda \in \Lambda$, 
\[
a + \alpha \cdot \be_{A^{+}} = (a+\alpha \cdot \be_{A^{+}})^{\halb}(\be_{A^{+}} - u_{\lambda})  (a+\alpha \cdot \be_{A^{+}})^{\halb} + (a+\alpha \cdot \be_{A^{+}})^{\halb} u_{\lambda}  (a+\alpha \cdot \be_{A^{+}})^{\halb},
\]
with the second summand being an element of $A$ dominated by $a+ \alpha \cdot \be_{A^{+}}$. This yields that
\[
b \perp  (a+\alpha \cdot \be_{A^{+}})^{\halb} u_{\lambda}  (a+\alpha \cdot \be_{A^{+}})^{\halb},
\]
hence 
\begin{equation}
\label{ww1}
\varphi^{(+)}(b)  \varphi^{(+)}(  (a+\alpha \cdot \be_{A^{+}})^{\halb} u_{\lambda}  (a+\alpha \cdot \be_{A^{+}})^{\halb}) = 0,
\end{equation}
since $\varphi^{(+)}$ agrees with the order zero map $\varphi$ on $A$. 

Furthermore, using continuity of $\varphi^{(+)}$ and the fact that $(u_{\lambda})_{\Lambda}$ is approximately central with respect to $A^{+}$, we check that 
\begin{eqnarray}
0 & \le & \mathrm{s.o.} \lim_{\lambda} \varphi^{(+)}((a+\alpha \cdot \be_{A^{+}})^{\halb} (\be_{A^{+}} - u_{\lambda})  (a+\alpha \cdot \be_{A^{+}})^{\halb})  \nonumber \\
& = & \mathrm{s.o.} \lim_{\lambda} \varphi^{(+)}( (\be_{A^{+}} - u_{\lambda})  (a+\alpha \cdot \be_{A^{+}})) \nonumber \\
& = & \mathrm{s.o.} \lim_{\lambda} \varphi^{(+)}(\alpha \cdot (\be_{A^{+}} - u_{\lambda})) \nonumber  \\
& = & \alpha \cdot (\varphi^{(+)}(\be_{A^{+}}) -  \mathrm{s.o.} \lim_{\lambda} \varphi ( u_{\lambda}) ) \nonumber \\
& \stackrel{\eqref{w3},\eqref{w4}}{=} & 0. \label{ww2}
\end{eqnarray}
We obtain
\begin{eqnarray*}
\lefteqn{\varphi^{(+)}(b) \varphi^{(+)}(a + \alpha \cdot \be_{A^{+}})} \\
 & = &  \varphi^{(+)}(b) \varphi^{(+)}((a + \alpha \cdot \be_{A^{+}})^{\halb} u_{\lambda} (a + \alpha \cdot \be_{A^{+}})^{\halb} ) \\
& & +  \varphi^{(+)}(b) \varphi^{(+)}((a + \alpha \cdot \be_{A^{+}})^{\halb} (\be_{A^{+}} - u_{\lambda}) (a + \alpha \cdot \be_{A^{+}})^{\halb} ) \\
& \stackrel{\eqref{ww1}}{=} &  \varphi^{(+)}(b) \varphi^{(+)}((a + \alpha \cdot \be_{A^{+}})^{\halb} (\be_{A^{+}} - u_{\lambda}) (a + \alpha \cdot \be_{A^{+}})^{\halb} ) \\
& \stackrel{\mathrm{s.o.}}{\to} & 0
\end{eqnarray*}
(where the last assertion follows from \eqref{ww2}), which implies that 
\[
\varphi^{(+)}(b) \varphi^{(+)}(a + \alpha \cdot \be_{A^{+}}) =0.
\]
Therefore, $\varphi^{(+)}$ has order zero.

To show that $\varphi^{(+)}$ is the unique c.p.c.\ order zero extension of $\varphi$ (mapping from $A^{+}$ to $C^{**}$), suppose $\psi:A^{+} \to C^{**}$ was another such extension, with $d:= \psi(\be_{A^{+}})$. Since $\psi$ is positive, it is clear that 
\[
d = \psi(\be_{A^{+}}) \ge \mathrm{s.o.} \lim_{\lambda} \psi(u_{\lambda}) = g.
\]
Now, suppose that $\|d-g\| >0$. Using that $\varphi (u_{\lambda})^{\frac{1}{n}} \to \be_{C^{**}}$ strongly as $\lambda \to \infty$ and $n \to \infty$, it is straightforward to show that there are $\eta>0$ and $\lambda \in \Lambda$ such that 
\[
\|(d-g) \varphi(u_{\lambda})(d-g)\| \ge \eta.
\]
Using functional calculus, one finds $u,w \in A_{+}$ of norm at most one such that 
\[
\|u-u_{\lambda}\| < \eta/2
\] 
and 
\[
wu = u.
\]
The latter implies that $\be_{A^{+}}-w$ and $u$ are orthogonal elements in $A^{+}$, whence 
\[
\psi(\be_{A^{+}}-w) \perp \psi(u) = \varphi(u) = \varphi^{(+)}(u) \perp  \varphi^{(+)}(\be_{A^{+}}-w).
\]
Combining these facts, we obtain
\begin{eqnarray*}
\eta & \le & \|(d-g) \varphi(u_{\lambda})(d-g)\| \\
& \le & \|(d-g) \varphi(u_{\lambda})\| \\
& \le & \|(\psi(\be_{A^{+}}) - \varphi^{(+)}(\be_{A^{+}})) \varphi(u)\| + \frac{\eta}{2} \\
& = &  \|(\psi(\be_{A^{+}} - (\be_{A^{+}}-w)) - \varphi^{(+)}(\be_{A^{+}}- (\be_{A^{+}}-w))) \varphi(u)\| + \frac{\eta}{2} \\
& = & \|(\psi(w) - \varphi^{(+)}(w)) \varphi(u) \| + \frac{\eta}{2} \\
& = & \|(\varphi(w) - \varphi(w)) \varphi(u)\| + \frac{\eta}{2} \\
& = &  \frac{\eta}{2}, \\
\end{eqnarray*}
a contradiction, so that $d = g$ and $\psi$ and $\varphi^{(+)}$ coincide.
\end{nproof}
\en

\begin{rems}
It will follow from (the proof of) the next theorem that the range of the map $\varphi^{(+)}$ of the preceding proposition in fact lies in $\Mh(C)$. 
\end{rems}

\bn
\label{main-result}
Our main result is the following structure theorem for c.p.\ order zero maps.

\begin{thms}
Let $A$ and $B$ be $C^{*}$-algebras and $\varphi:A \to B$ a c.p.\ order zero map. Set $C:= C^{*}(\varphi(A)) \subset B$. 

Then, there are a positive element $h \in \Mh(C) \cap C'$ with $\| h \| = \|\varphi \|$ and a $*$-homomorphism
\[
\pi_{\varphi}:A \to \Mh(C) \cap \{h\}'
\]
such that
\[
\pi_{\varphi}(a)h = \varphi(a) \mbox{ for } a \in A.
\]
If $A$ is unital, then $h = \varphi(\be_{A}) \in C$.
\end{thms}

\begin{nproof}
By rescaling $\varphi$ if necessary, we may clearly assume that $\varphi$ is contractive. Let us first assume that $A$ is unital, and set 
\begin{equation}
\label{w9}
h:= \varphi(\be_{A}) \in C.
\end{equation}
We may further assume that $C$ acts nondegenerate on its universal Hilbert space $\mathcal{H}$. 

By \cite[Theorem~2.3(i)]{Wol:disjointness} (cf.\ \ref{Wolffs-result}), we have $h \in \mathcal{Z}(C)$, since $A$ is unital and $\varphi$ is disjointness preserving. Moreover, one checks that $h$ is a strictly positive element of $C$, and since $C \subset \mathcal{B}(\mathcal{H})$ is nondegenerate, this implies that the support projection of $h$ is $\be_{\mathcal{H}}$. On the other hand, the support projection of $h$ can be expressed as 
\begin{equation}
\label{w1}
\mathrm{s.o.}\lim_{n \to \infty} (h+ \frac{1}{n} \cdot \be_{\mathcal{H}})^{-1} h = \be_{\mathcal{H}} = \be_{C^{**}}.
\end{equation}
We now define a map
\[
\pi_{\varphi} : A \to C^{**} \subset \mathcal{B}(\mathcal{H})
\]
by
\begin{equation}
\label{w8}
\pi_{\varphi}(a):= \mathrm{s.o.}\lim_{n \to \infty} (h+ \frac{1}{n} \cdot \be_{\mathcal{H}})^{-1} \varphi(a).
\end{equation}
Existence of the limit can be checked on positive elements, since then the sequence $(h+\frac{1}{n} \cdot \be_{\mathcal{H}})^{-1}\varphi(a)$ is monotone increasing. Since $\pi_{\varphi}$ is a strong limit of c.p.\ maps, it is c.p.\ itself. Since $h$ commutes with $\varphi(A)$, one checks that $\pi_{\varphi}$ again has order zero. Moreover, 
\begin{equation}
\label{w10}
\pi_{\varphi}(\be_{A}) = \be_{\mathcal{H}}
\end{equation} 
by \eqref{w1} and \eqref{w8}, so $\pi_{\varphi}$ is unital. 

Now by \cite[Lemma~3.3]{Wol:disjointness}, $\pi_{\varphi}$ is a Jordan $*$-homomorphism, so 
\begin{equation}
\label{w6}
\pi_{\varphi}(a^{2}) = \pi_{\varphi}(a)^{2} \mbox{ for } a \in A_{+}.
\end{equation}
On the other hand, since $\pi_{\varphi}$ is u.c.p., by Stinespring's Theorem we may assume that there is a unital $C^{*}$-algebra $D$ containing $C^{**}$ and a $*$-homomorphism 
\[
\varrho: A \to D
\]
such that 
\begin{equation}
\label{w5}
\pi_{\varphi}(a) = \be_{C^{**}} \varrho(a) \be_{C^{**}} \mbox{ for } a \in A.
\end{equation}
We now compute 
\begin{eqnarray}
\lefteqn{\| \be_{C^{**}} \varrho(a) - \be_{C^{**}} \varrho(a) \be_{C^{**}}\|^{2}} \nonumber \\
& = & \| \be_{C^{**}} \varrho(a) (\be_{D} - \be_{C^{**}}) \varrho(a) \be_{C^{**}} \| \nonumber \\
& \stackrel{\eqref{w5}}{=} & \| \pi_{\varphi}(a^{2}) - \pi_{\varphi}(a)^{2} \| \nonumber \\
& \stackrel{\eqref{w6}}{=} & 0 \label{w7}
\end{eqnarray}
for $a \in A_{+}$, whence 
\begin{eqnarray*}
\pi_{\varphi}(ab) & \stackrel{\eqref{w5}}{=} & \be_{C^{**}} \varrho(a) \varrho(b) \be_{C^{**}} \\
& \stackrel{\eqref{w7}}{=} &  \be_{C^{**}} \varrho(a) \be_{C^{**}} \varrho(b) \be_{C^{**}} \\
& \stackrel{\eqref{w5}}{=} & \pi_{\varphi}(a) \pi_{\varphi}(b)
\end{eqnarray*}
for $a,b \in A_{+}$. By linearity of $\pi_{\varphi}$ it follows that $\pi_{\varphi}$ is multiplicative, hence a $*$-homomorphism.

Next, we check that for $a \in A$,
\begin{eqnarray*}
\varphi(a) - \pi_{\varphi}(a) h & \stackrel{\eqref{w8}}{=} & \varphi(a) - \varphi(a) \mathrm{s.o.} \lim_{n \to \infty} (h + \frac{1}{n} \cdot \be_{\mathcal{H}})^{-1} h \\
& \stackrel{\eqref{w9}}{=} & \varphi(a) - \varphi(a)   \mathrm{s.o.} \lim_{n \to \infty} (h + \frac{1}{n} \cdot \be_{\mathcal{H}})^{-1} \varphi(\be_{A}) \\
& \stackrel{\eqref{w8}}{=} & \varphi(a) - \varphi(a) \pi_{\varphi}(\be_{A}) \\
& \stackrel{\eqref{w10}}{=} & 0,
\end{eqnarray*}
so
\[
\varphi(a)= \pi_{\varphi}(a)h = h \pi_{\varphi}(a)
\]
for all $ a \in A$, and $\pi_{\varphi}(A) \subset \{h\}'$.

Finally, we have for $a,b \in A$ 
\[
\pi_{\varphi}(a) \varphi(b) = \pi_{\varphi}(a) \pi_{\varphi}(b) h = \pi_{\varphi}(ab) h = \varphi(ab) \in C,
\]
and similarly $\varphi(a) \pi_{\varphi}(b)=\varphi(ab)$, 
from which one easily deduces that
\[
\pi_{\varphi}(A) \subset \Mh(C).
\]
We have now verified the lemma in the case where $A$ is unital. In the nonunital case, we may use Proposition~\ref{order-zero-unitization} to extend $\varphi$ to a c.p.c.\ order zero map $\varphi^{(+)}:A^{+} \to C^{**}$. By the first part of the proof there is a $*$-homomorphism $\pi_{\varphi^{(+)}}: A^{+} \to C^{**}$ such that $\varphi^{(+)}(a)  = \pi_{\varphi^{(+)}}(a) g = g \pi_{\varphi^{(+)}}(a)$ for all $a \in A^{+}$, where $g := \varphi^{(+)}(\be_{A^{+}})$. Now if $b \in A_{+}$, we have 
\begin{eqnarray*}
g \varphi(b) & = & g \varphi^{(+)}(b) \\
& = & g \pi_{\varphi^{(+)}}(b) g \\
& = & g \pi_{\varphi^{(+)}}(b^{\halb}) \pi_{\varphi^{(+)}}(b^{\halb}) g \\
& = & \varphi^{(+)}(b^{\halb}) \varphi^{(+)}(b^{\halb}) \\
& = & \varphi(b^{\halb})^{2} \in C,
\end{eqnarray*}
which, by linearity, yields $g \varphi(b) \in C$ for any $b \in A$. From here it is straightforward to conclude that  $g \in \Mh(C)$, whence the images of $\varphi^{(+)}$ and $\pi_{\varphi^{(+)}}$ in fact both live in $\Mh(C)$ by the first part of the proof. The $*$-homomorphism $\pi_{\varphi }:A \to \Mh(C)$ will then just be the restriction of $\pi_{\varphi^{(+)}}:A^{+} \to \Mh(C)$  to $A$. 
\end{nproof}
\en

\section{Some consequences}
\label{consequences}

In this final section we derive some corollaries from Theorem~\ref{main-result}.

\bn
\begin{cors}
Let $A$ and $B$ be $C^{*}$-algebras, and $\varphi:A \to B$ a c.p.c.\ order zero map. Then, the map given by $\varrho_{\varphi}(\id_{(0,1]} \otimes a) := \varphi(a)$ (for $a \in A$) induces a $*$-homomorphism $\varrho_{\varphi}: \Ch_{0}((0,1]) \otimes A \to B$.

Conversely, any $*$-homomorphism $\varrho: \Ch_{0}((0,1]) \otimes A \to B$ induces a c.p.c.\  order zero zero map $\varphi_{\varrho}: A \to B$ via $\varphi_{\varrho}(a):= \varrho(\id_{(0,1]} \otimes a)$.

These mutual assignments yield a canonical bijection between the spaces of c.p.c.\ order zero maps from $A$ to $B$ and $*$-homomorphisms from $\Ch_{0}((0,1]) \otimes A$ to $B$.
\end{cors}

\begin{nproof}
It is well known that $\Ch_{0}((0,1])$ is canonically isomorphic to the universal $C^{*}$-algebra generated by a positive contraction, identifying $\id_{(0,1]}$ with the universal generator. 

Now if $\varphi: A \to B$ is c.p.c.\ order zero, obtain $C$, $h$ and $\pi_{\varphi}$ from Theorem~\ref{main-result}. There is a $*$-homomorphism
\[
\bar{\varrho}: \Ch_{0}((0,1]) \to \Mh(C)
\]
induced by
\[
\bar{\varrho}(\id_{(0,1]}):= h;
\]
since $h \in \pi_{\varphi}(A)'$, $\bar{\varrho}$ and $\pi_{\varphi}$ yield a $*$-homomorphism 
\[
\varrho_{\varphi}: \Ch_{0}((0,1]) \otimes A \to \Mh(C)
\]
satisfying
\[
\varrho_{\varphi} (\id_{(0,1]} \otimes a) = h \pi_{\varphi}(a) = \varphi(a) \in C
\]
for $a \in A$. Since $\Ch_{0}((0,1]) \otimes A$ is generated by $\id_{(0,1]} \otimes A$ as a $C^{*}$-algebra, we see that in fact the image of $\varrho_{\varphi}$ lies in $C \subset B$. 

Conversely, if 
\[
\varrho: \Ch_{0}((0,1]) \otimes A \to B
\]  
is a $*$-homomorphism, then
\[
\varphi_{\varrho}(\, . \,):= \varrho \circ (\id_{(0,1]} \otimes \, . \,)
\]
clearly has order zero. 

That the assignments $\varphi \mapsto \varrho_{\varphi}$ and $\varrho \mapsto \varphi_{\varrho}$ are mutual inverses is straightforward to check.
\end{nproof}
\en

\bn
\label{order-zero-functional-calculus}
As in \cite{Winter:dr-Z-stable}, Theorem~\ref{main-result} allows us to define a positive functional calculus of c.p.c.\ order zero maps.

\begin{cors}
Let $\varphi:A \to B$ be a c.p.c.\ order zero map, and let $f \in \Ch_{0}((0,1])$ be a positive function. Let $C$, $h$ and $\pi_{\varphi}$ be as in Theorem~\ref{main-result}. Then, the map 
\[
f(\varphi): A \to C \subset B,
\]
given by
\[
f(\varphi)(a):= f(h) \pi_{\varphi}(a) \mbox{ for } a \in A,
\]
is a well-defined c.p.\ order zero map. If $f$ has norm at most one, then $f(\varphi)$ is also contractive.
\end{cors}

\begin{nproof}
Since $[h,\pi_{\varphi}(A)]=0$, we also have $[f(h),\pi_{\varphi}(A)]=0$, which implies that $f(\varphi)$ indeed is a c.p.\ map. Using that $h \pi_{\varphi}(a) \in C$ for any $a \in A$, it is straightforward to conclude that $f(h) \pi_{\varphi}(a) \in C$ for any $a \in A$. The last statement is obvious.
\end{nproof}
\en

\bn
\label{tensor}
\begin{cors}
Let $A$, $B$, $C$ and $D$ be $C^{*}$-algebras and $\varphi:A \to B$ and $\psi:C \to D$ c.p.c.\ order zero maps. 

Then, the induced c.p.c.\ map 
\[
\varphi \otimes_{\mu} \psi: A \otimes_{\mu} C \to B \otimes_{\mu} D
\]
has order zero, when $\otimes_{\mu}$ denotes the minimal or the maximal tensor product. In particular, for any $k \in \mathbb{N}$ the amplification
\[
\varphi^{(k)}: M_{k}(A) \to M_{k}(B)
\]
has order zero.
\end{cors}

\begin{nproof}
Set \[
\bar{B}:= C^{*}(\varphi(A)) \subset B \mbox{ and } \bar{D}:= C^{*}(\psi(C)) \subset D,
\]
and employ Theorem~\ref{main-result} to obtain $*$-homomorphisms
\[
\pi_{\varphi}:A \to \Mh(\bar{B}) \mbox{ and } \pi_{\psi}: C \to \Mh(\bar{D})
\]
and positive elements
\[
h_{\varphi} \in \Mh(\bar{B}) \mbox{ and } h_{\psi} \in \Mh(\bar{D}),
\]
so that
\[
\varphi(a)= \pi_{\varphi}(a) h_{\varphi } = h_{\varphi} \pi_{\varphi}(a) \mbox{ and } \psi(a)= \pi_{\psi}(c) h_{\psi } = h_{\psi} \pi_{\psi}(c)
\]
for $a \in A$ and $c \in C$. 

Let us consider the maximal tensor product first. Fix a faithful representation 
\[
\iota: B \otimes_{\mathrm{max}} D \hookrightarrow \Bh(\Hh);
\]
we have 
\[
\iota_{\mathrm{max}}: \bar{B} \otimes_{\mathrm{max}} \bar{D} \to \bar{B} \otimes_{\nu} \bar{D} \subset B \otimes_{\mathrm{max}} D \subset \Bh(\Hh)
\]
for some $C^{*}$-norm $\nu$ on $\bar{B} \odot \bar{D}$ ($\iota_{\mathrm{max}}$ is not necessarily injective). The representation of $\bar{B} \odot \bar{D}$ on $\Hh$ yields representations of $\bar{B}$ and $\bar{D}$ on $\Hh$ with commuting images (cf.\ \cite[Theorem~II.9.2.1]{Bla:encyc}), and one checks that the induced representations of $\Mh(\bar{B})$ and $\Mh(\bar{D})$ on $\Hh$ also commute, and live in $\Mh(\bar{B} \otimes_{\nu} \bar{D})$. We then obtain a $*$-homomorphism
\[
\bar{\iota}_{\mathrm{max}}: \Mh(\bar{B}) \otimes_{\mathrm{max}} \Mh(\bar{D}) \to \Mh(\bar{B} \otimes_{\nu} \bar{D}) \subset \Bh(\Hh)
\]
extending $\iota_{\mathrm{max}}$, so that we may define a $*$-homomorphism
\[
\pi_{\mathrm{max}}: A \otimes_{\mathrm{max}} C \stackrel{\pi_{\varphi} \otimes_{\mathrm{max}} \pi_{\psi}}{\longrightarrow} \Mh(\bar{B}) \otimes_{\mathrm{max}} \Mh(\bar{D}) \stackrel{\bar{\iota}_{\mathrm{max}}}{\longrightarrow} \Mh(\bar{B} \otimes_{\nu} \bar{D}) \subset \Bh(\Hh)
\]
and a positive element
\[
h_{\mathrm{max}}:= \bar{\iota}_{\mathrm{max}} \circ (h_{\varphi} \otimes h_{\psi}) \in \Mh(\bar{B} \otimes_{\nu} \bar{D}) \subset \Bh(\Hh).
\]
It is straightforward to verify that
\[
[h_{\mathrm{max}}, \pi_{\mathrm{max}}(A \otimes_{\mathrm{max}} C)] = 0
\]
and that
\[
\varphi \otimes_{\mathrm{max}} \psi : A \otimes_{\mathrm{max}} C \to \bar{B} \otimes_{\mathrm{max}} \bar{D} \to \bar{B} \otimes_{\nu} \bar{D} \subset B \otimes_{\mathrm{max}} D
\]
satisfies
\[
\varphi \otimes_{\mathrm{max}} \psi(\, . \,) = \pi_{\mathrm{max}}(\, .\,) h_{\mathrm{max}},
\]
so $\varphi \otimes_{\mathrm{max}} \psi$ indeed has order zero. 

The minimal tensor product is handled similarly; only now we have to consider faithful representations $\iota_{\bar{B}}: \bar{B} \to \Bh(\Hh_{\bar{B}})$ and $\iota_{\bar{D}}: \bar{D} \to \Bh(\Hh_{\bar{D}})$ which induce a (faithful) representation $\bar{\iota}_{\mathrm{min}}: \Mh(\bar{B}) \otimes_{\mathrm{min}} \Mh(\bar{D}) \to \Bh(\Hh_{\bar{B}} \otimes \Hh_{\bar{D}})$. As above, we set
\[
\pi_{\mathrm{min}}:= \bar{\iota}_{\mathrm{min}} \circ ( \pi_{\varphi} \otimes_{\mathrm{\min}} \pi_{\psi}): A \otimes_{\mathrm{min}} C \to \Bh(\Hh_{\bar{B}} \otimes \Hh_{\bar{D}})
\]
and 
\[
h_{\mathrm{\min}}:= \bar{\iota}_{\mathrm{min}} \circ (h_{\varphi} \otimes h_{\psi}),
\]
and check that 
\[
[h_{\mathrm{min}}, \pi_{\mathrm{min}}(A \otimes_{\mathrm{min}} C)] = 0
\]
and
\[
\varphi \otimes_{\mathrm{min}} \psi(\, . \,) = \pi_{\mathrm{min}}(\, . \,) h_{\mathrm{min}}.
\]
\end{nproof}
\en

\bn
\begin{cors}
Let $A$ and $B$ be $C^{*}$-algebras, $\varphi:A \to B$ a c.p.c.\ order zero map, and $\tau$ a positive tracial functional on $B$. 

Then, the composition $\tau \circ \varphi$ is a positive tracial functional. 

The statement also holds when replacing the term `positive tracial functional' with `2-quasitrace'.
\end{cors}

\begin{nproof}
If $\tau$ is a positive tracial functional, we only need to check that $\tau \varphi$ satisfies the trace property. But, for $a,b \in A$, 
\begin{eqnarray*}
\tau \varphi(ab) & = & \tau(\varphi^{\halb}(a) \varphi^{\halb}(b) )  \\
& = &  \tau(\varphi^{\halb}(b) \varphi^{\halb}(a) )  \\
& = & \tau \varphi(ba).
\end{eqnarray*}
Here, we have used the trace property of $\tau$ as well as \ref{order-zero-functional-calculus}. 

If $\tau$ is only a 2-quasitrace, we also have to check two other things: For once, that $\tau \varphi$ extends to $M_{2}(A)$ -- but this is obvious as $\varphi$ extends to a c.p.c.\ order zero map by Corollary~\ref{tensor}. Second, we need to check that $\tau \varphi$ is additive on commuting elements. However, if $a,b \in A$ satisfy $[a,b]=0$, then $[\varphi(a),\varphi(b)]= h^{2}\pi_{\varphi}([a,b])=0$ and
\[
\tau\varphi(a+b) = \tau(\varphi(a)+\varphi(b)) = \tau\varphi(a) + \tau\varphi(b),
\]
since $\varphi$ is linear and $\tau$ is a 2-quasitrace.
\end{nproof}
\en

\bn
The next result is one of our main motivations for studying order zero maps in the abstract; it says that order zero maps induce maps at the level of Cuntz semigroups, since they respect the Cuntz relation; see \cite{Cuntz:dimension}, \cite{Ror:uhfII} and  \cite{Ror:Z-absorbing} for an introduction to Cuntz subequivalence and the Cuntz semigroup.

\begin{cors}
Let $A$ and $B$ be $C^{*}$-algebras and $\varphi:A \to B$ a c.p.c.\ order zero map. 

Then, $\varphi$ induces a morphism of ordered semigroups 
\[
W(\varphi): W(A) \to W(B)
\]
between the Cuntz semigroups via 
\[
W(\varphi)(\langle a \rangle) = \langle \varphi^{(k)}(a) \rangle \mbox{ if } a \in M_{k}(A)_{+}.
\]
\end{cors}

\begin{nproof}
Let $a,b \in M_{k}(A)_{+}$ for some $k \in \mathbb{N}$ (it clearly suffices to consider the same $k$ for $a$ and $b$) satisfying $a \precsim b$. By definition of Cuntz subequivalence (cf.\ \cite{Cuntz:dimension}), this means that there is a sequence $(x_{n})_{\mathbb{N}} \subset M_{k}(A)$ such that 
\[
a = \lim_{n \to \infty} x_{n}^{*} b x_{n}.
\]
Let
\[
\varphi^{(k)}: M_{k}(A) \to M_{k}(B)
\]
denote the amplification of $\varphi$; note that $\varphi^{(k)}$ has order zero by Corollary~\ref{tensor}. Let 
\[
C:= C^{*}(\varphi^{(k)}(M_{k}(A))) (\cong M_{k}(C^{*}(\varphi(A))))
\]
and obtain from Theorem~\ref{main-result} a $*$-homomorphism 
\[
\pi_{\varphi^{(k)}}: M_k(A) \to \Mh(C)
\]
and
\[
h \in \Mh(C)_{+}
\]
commuting with $C$. For $n \in \mathbb{N}$, define 
\[
\tilde{x}_{n}:= h^{\frac{1}{n}} \pi_{\varphi^{(k)}}(x_{n}) = (\varphi^{(k)})^{\frac{1}{n}}(x_{n})
\]
using \ref{order-zero-functional-calculus}; we then have  
\[
\tilde{x}_{n} \in C,
\]
and
\begin{eqnarray*}
\varphi^{(k)}(a) & = & \lim_{n \to \infty} \varphi^{(k)} (x_{n}^{*} b x_{n}) \\
& = & \lim_{n \to \infty} \pi_{\varphi^{(k)}}(x_{n}^{*})h^{\halb} \pi_{\varphi^{(k)}}(b) h^{\halb}\pi_{\varphi^{(k)}}(x_{n}) \\
& = & \lim_{n \to \infty} \tilde{x}_{n}^{*} \varphi^{(k)}(b) \tilde{x}_{n},
\end{eqnarray*}
where we have used that 
\[
\lim_{n \to \infty} h^{\halb} h^{\frac{1}{n}} = h^{\halb}.
\]
It follows that
\[
\langle \varphi^{(k)}(a) \rangle \le \langle \varphi^{(k)}(b) \rangle,
\]
so that 
\[
W(\varphi)(\langle a \rangle) \le W(\varphi)(\langle b \rangle).
\]
The argument also shows that, if $a \sim b$, then $\varphi^{(k)}(a) \sim \varphi^{(k)}(b)$, so that $W(\varphi)$ indeed is well-defined and respects the order. 
Moreover, if $a,b \in M_{k}(A)$ are orthogonal, then so are $\varphi^{(k)}(a),\varphi^{(k)}(b) \in M_{k}(B)$, whence
\[
\varphi^{(k)}(a \oplus b) = \varphi^{(k)}(a) \oplus \varphi^{(k)}(b)
\]
and $W(\varphi)$ is a semigroup morphism.
\end{nproof}
\en

\bn
We remark in closing that if $\Mh$ and $\Nh$ are von Neumann algebras and $\varphi: \Mh \to \Nh$ is a c.p.c .\ order zero map then the proof of \ref{main-result} shows that $\varphi = h \pi_{\varphi}$
where $h \in \Nh_+$ commutes with the range of $\varphi$ and $\pi_{\varphi}$ is a $*$-homomorphism which is normal if $\varphi$ is normal.
Moreover, if $\varphi: A \to B$ is any c.p.c.\ order zero map between $C^*$-algebras then so is its bitransposed $\varphi^{**}: A^{**} \to B^{**}$.
This follows for instance by bitransposing the factorization $\varphi= (h^{1/2} \cdot h^{1/2}) \circ  \pi_{\varphi}$ and using $\Mh(C)^{**}= C^{**} \oplus ( \Mh(C)/C)^{**}$ (if $A$ is nonunital, one has to also use Proposition~\ref{order-zero-unitization}, and in particular \eqref{w4}).
\en





%

%

%

%


\providecommand{\bysame}{\leavevmode\hbox to3em{\hrulefill}\thinspace}
\providecommand{\MR}{\relax\ifhmode\unskip\space\fi MR }
\providecommand{\MRhref}[2]{%
  \href{http://www.ams.org/mathscinet-getitem?mr=#1}{#2}
}
\providecommand{\href}[2]{#2}

\end{document}